\documentclass{elsart}

\begin{document}

\begin{frontmatter}



\title{Definability as
Hypercomputational Effect}


\author{S.\ Barry Cooper\thanksref{label1}}
\thanks[label1]{The author wishes to thank Istv\'an N\'emeti 
for a number of helpful comments on a first draft of this paper, 
which have  led to a number of improvements. \endgraf
Preparation of this article partially 
supported by EPSRC Research Grant EP/C001389/1 
{\it Extensions of Embeddings in the Local Turing Universe\/}, and by 
NSFC Grand International Joint Project Grant 
No.\ 60310213 {\it New Directions in
Theory and Applications of Models of Computation.}}

\address{School of Mathematics, University of Leeds\\
Leeds LS2 9JT, U.K.}

\begin{abstract}The classical simulation of physical 
processes using standard models of computation is 
fraught with problems. On the other hand, attempts at 
modelling real-world computation with the aim 
of isolating its hypercomputational content have 
struggled to convince. We argue that a better basic 
understanding can be achieved through 
computability theoretic deconstruction of 
those physical phenomena most resistant 
to classical simulation.  
From this we may be able to better assess whether the  
hypercomputational enterprise is proleptic computer 
science, or of mainly philosophical interest. \end{abstract}

\end{frontmatter}

\section{Introduction}\label{intro}

In the face of new challenges to familiar computational paradigms, 
the role of classical computability theory, based on Turing machines and the 
Church-Turing Thesis, has appeared to be a largely reactive one. 
Reflective of the growing body of literature openly questioning 
the relevance of classical models of computation is the 
following quotation from Goldin and Wegner's \cite{GW} paper from 
CiE 2005:

{\vskip 0mm\noindent\rightskip=14pt \leftskip=14pt \hskip 0mm {
``It is time to recognise that today's computing applications, 
such as web services, intelligent agents, operating systems, and 
graphical user interfaces, cannot be modeled by 
Turing machines; alternative models are needed."\vskip
0mm}}

This echoes and extends van Leeuwen and  
Wiedermann's \cite{LW}  
 observation\footnote{See also the related paper by the same authors
on {\it Relativistic computers and non-uniform
complexity theory.} In: Calude et al (eds.) UMC 2002. Lecture Notes in
Computer Science Vol. 2509, Springer-Verlag, Berlin, 2002. pp. 287--299.}
that ``the classical Turing
paradigm may no longer be fully appropriate  
to capture all features of present-day
computing.''

In contrast to this, one has Martin Davis 
\cite{davis04} confidently asserting that: 

{\vskip 0mm\noindent\rightskip=14pt \leftskip=14pt \hskip 0mm {
``The great success of modern computers as all-purpose 
algorithm-executing engines embodying Turing's 
universal computer in physical form, makes it extremely 
plausible that the abstract theory of computability gives the 
correct answer to the question `What is a computation?', and, 
by itself, makes the existence of any more general form 
of computation extremely doubtful."\vskip 0mm}}

This has become a hard position to 
maintain. If only, as Philip Welch \cite{welch} points out, because  
\textit{logical\/}  
proofs of the impossibility of hypercomputation ``may be akin to 
proofs that fairies are logically impossible: damn hard to be convincing." 

In 1998 Jack Copeland \cite{copeland98} claimed to have rediscovered in
Turing's  1939 paper \cite{turing39}, based 
on his 1938 thesis at Princeton,  
a previously unsuspected hypercomputational agenda, 
based on the familiar oracle Turing machine --- an exegesis which seemed to  
throw the embattled classical computability theorist 
a lifeline. But Davis \cite{davis04} was having none of it: 

{\vskip 0mm\noindent\rightskip=14pt \leftskip=14pt \hskip 0mm {
``It is perfectly plain in the context of Turing's dissertation, 
that O-machines were introduced simply to solve a specific 
technical problem about definability of sets of natural 
numbers. There is not the faintest hint that Turing was making 
a proposal about a machine to be built. \dots 
It makes no sense to imagine that he was thinking about actual 
machines to compute the uncomputable."\vskip 0mm}}

But just say we \textit{do\/} have misgivings about 
universality of the universal Turing machine --- surely it makes 
sense to look more closely at 
problematic computational situations, using all the 
technical sophistication which computability 
theory puts at our  disposal. May that not lead to 
computational models with sufficient explanatory power 
to throw light on current controversies?  Whatever 
our interpretation of Turing's real views on 
mechanical computability (and in \cite{cooper05}  
I have indicated that these were not as simple as sometimes 
thought), Copeland has effectively spotlighted 
the potential for mathematically analysing the computationally 
complex in terms of its algorithmic content. As we argued 
in \cite{cooper98}: 

{\vskip 0mm\noindent\rightskip=14pt \leftskip=14pt \hskip 0mm {
``If one abstracts from the Universe its information 
content, structured via the basic \dots 
 fundamental laws of nature, one 
obtains a particular \dots   
manifestation of the Turing universe \dots , within which vague questions 
attain a precise analogue of quite plausible
validity."\vskip 0mm}}

What is needed is for the same level of deconstructive 
analysis Turing applied to the 
intuitive concept of computation to be applied 
to  hypercomputation.  

This article is motivated 
by two important questions. One is the computability-theoretic 
question concerning the nature of the relationship between 
computation  and 
(hyper)computational effects. The other 
is the more down-to-earth question of whether hypercomputation, 
if it actually occurs, can be successfully harnessed to 
human computational needs. 
This is how Martin Davis \cite{davis} 
puts it in relation to Etesi and Nemeti's discussion \cite{EN} 
of relativistic behaviour near black holes: 

{\vskip 0mm\noindent\rightskip=14pt \leftskip=14pt \hskip 0mm {
``Of course, even assuming that all this really does correspond 
to the actual universe in which we live, there is still the 
question of whether an actual device to take advantage of this 
phenomenon is possible."\footnote{Davis adding ``But the theoretical question
is certainly of
 interest."}\vskip 0mm}}

Or, more negatively (Davis \cite{davis04}): 

{\vskip 0mm\noindent\rightskip=14pt \leftskip=14pt \hskip 0mm {
``In any case, a useable physical representation of an 
uncomputable function, would require a revolutionary new 
physical theory, and one that it would be impossible to 
verify because of the inherent limitations of physical 
measurement."\vskip 0mm}}

Is Davis right?\footnote{Davis' particular doubts here are
addressed on a number of  levels in the 
Istv\'an N\'emeti and Gyula D\'avid article \cite{ND} in this 
volume. }
 Of course, in principle, he might well be. 
What may change that is a better understanding of the computational 
character of physical processes. In which case, the answer to the second
question will give us a better idea of whether the 
hypercomputational enterprise is proleptic computer 
science, or philosophy. 

The approach of  
 N\'emeti and  D\'avid \cite{ND} is a rather different one from 
that taken here. Theirs is based on what they describe as `a major 
paradigm shift in our physical world-view as well as our cosmological 
one'. They seek to give {\it scientific\/} substance to the 
`fairies' Welch refers to: an enterprise which, 
as was argued in \cite{CoOd03},  is `damn hard to be convincing', 
  but which we would 
expect to provide an important element in the establishment of 
a  
larger paradigm shift.

\section{Neural mapping and  real-world hypercomputation}\label{neural}

If one is looking for a context which is both widely 
suspected of transcending the standard Turing model, 
and of whose inner workings we have a high level of detailed 
knowledge, we need look no further than the human brain. 
Turing himself  \cite{turing54} talks about  ``the inadequacy of `reason'
unsupported  by common 
sense", and in his 1939 paper \cite{turing39} says:

{\rightskip=14pt \leftskip=14pt \hskip 0mm {
``Mathematical  reasoning may be regarded \dots as the exercise 
of a combination of \dots {\sl intuition\/} and 
{\sl ingenuity\/}. \dots In pre-G\"odel times it was thought by some that 
all the intuitive judgements of mathematics could be replaced by 
a finite number of \dots rules. The necessity for intuition would then be 
entirely eliminated. 
In our discussions, however, we have gone to the opposite 
extreme and eliminated not intuition but ingenuity"}\vskip 0mm}
 
Teuscher \cite{teuscher} identifies the gap to be filled by 
suitably extended models:

{\rightskip=14pt \leftskip=14pt \hskip 0mm {
``People have started thinking about the possibility that 
simulating the mind in silicon might be impossible --- 
or at least impossible using today's methods. Should we 
first forget about computers and look closer at the 
gray stuff in the brain, since the actual knowledge of the 
brain is severely fragmented and many details are still not at all 
understood?"}\vskip 0mm}

Part 
of the brain's potential for enrichment of our modelling of 
the computationally complex lies in the way it 
seems to successfully deal with the sort of imaging 
of the real world we would dearly like our 
computing machines to perform. More important, 
the brain shows the capacity to perform re-presentations 
of mental imaging to enable recursive 
development of complex conceptual structures. At the same time, 
new techniques for relating structural and functional 
features of the brain, for example, using positron emission 
scan (PET) or a functional magnetic resonance imaging scan (fMRI), 
bring us much closer to  obtaining useful models. 

Connectionist models of computation 
based on the workings of the human brain have come a long way since 
Turing's \cite{turing48} 
discussion of `unorganised  machines'   (see    
 Jack Copeland and Diane Proudfoot's article \cite{CP96}  
\textit{On Alan Turing's Anticipation of Connectionism}), and 
McCulloch and Pitts' seminal paper \cite{MP} on neural nets. 
But despite the growth of computational neuroscience 
as
an active research area, 
putting together ingredients from both artificial 
neural networks and neurophysiology, 
something does seem to be missing. 
This leads Rodney Brooks 
\cite{brooks} to allude to the fact  that 
``neither AI nor Alife has produced artifacts that could be confused with a 
living organism for more than an instant." 
Again:  
 ``\dots neural networks 
alone cannot do the job" opines Steven Pinker, 
going on to describe \cite[p.124]{pinker} ``a kind of mental 
fecundity called recursion":

{\rightskip=14pt \leftskip=14pt \hskip 0mm {
``We humans can take an entire proposition and give it a 
role in some larger proposition. Then we can take the larger 
proposition and embed it in a still-larger one. 
Not only did the baby eat the slug, but the father saw the baby eat
the slug, and I wonder whether the father 
saw the baby eat the slug,
the father knows that I wonder whether he 
saw the baby eat the slug, and I can guess that the father knows that 
I wonder whether he saw the baby eat the slug, and so on."}\vskip 0mm}

We are good at 
devising computational models capable of imaging, and of 
going some way to emulate how the brain comes up 
with neural patterns representing quite complex 
formations. But the mechanisms the brain uses to represent 
such patterns and relate them in complex ways is more elusive. 
What makes the sort of recursion Stephen Pinker 
has in mind so difficult to get to grips 
with at the structural level, is that it seems wound up 
with the puzzle of consciousness and its relationship 
to emotions  and   
 feelings. Antonio Damasio  \cite[p.169]{damasio}  has a nice description
of  the hierarchical development 
of a particular instance of consciousness within the 
brain (or, rather, `organism'), interacting with some external 
object: 

{\rightskip=14pt \leftskip=14pt \hskip 0mm {
``\dots both organism and object are mapped as 
neural patterns, in first-order maps; all of these 
neural patterns can become images. \dots 
The sensorimotor maps pertaining to the 
object cause changes in the maps pertaining to the 
organism. \dots [These] changes \dots can be 
re-represented in yet other maps (second-order maps) 
which thus represent the relationship of object 
and organism. \dots The neural patterns 
transiently formed in second-order maps can become 
mental images, no less so than the neural 
patterns in first-order maps."}\vskip 0mm}

Notice that what is envisaged is the re-representation 
of neural patterns formed across some region of the brain, 
in such a way that they can have a computational relevance 
in forming new patterns. This is where the 
clear demarcation between computation and computational 
effect  becomes blurred. The key conception is of 
computational loops incorporating these `second-order' 
aspects of the computation itself. 
Building on this 
one can derive a plausible schematic picture of 
of the global workings of the brain.

Considering how complex a structure the human brain is, 
it is surprising one does not find more features needing 
reflecting in any basic computational model based on it. 
However, a thorough trawl through the literature, and one's own 
experiences, 
fails to bring to light anything that might be 
held up as  computational principle transcending 
in a fundamental way what we have already identified. 
The key ingredients we expect in a model are imaging, 
parallelism, interconnectivity, and a counterpart to the  
second-order recursions pointed to above. 

Mathematically, 
the imaging appears to be dependent on the  
parallelism and  interconnectivity. This is what 
connectionist models are strong on.   
The recursions are not so easy to model, though. 
Looked at logically, 
one has representations of complex patternings of neural 
events underlying which there is no clear local mechanism, 
but for which one would expect a description in terms of 
the  structures pertaining. Looked at physically, 
such descriptions appear to emerge, and be associated with 
 (but not exclusively) the sort of non-linear mathematics governing the 
emergence of new relations from chaotic environments. 
This leads us to turn the picture of re-representations 
of mental imaging as a describable mapping on its head, 
and think (see \cite{cooper.emerg}) in terms of descriptions 
in terms of a structure \textit{defining}, and hence 
determining, the mental re-representations. 

Looking at this more closely, what seems to be 
happening is that the brain stores away not just the 
image, but a route to accessing that image as a whole. 
This is what people who specialise in memorising 
very long numbers seem to display --- rather than attempting to go directly 
into the detailed memory of a given number, they use 
simple representational tricks to call the entire 
number up. Here is how Damasio summarises the process (and 
the quotation from \cite[p.170]{damasio}  
  is worth giving in full): 

{\rightskip=14pt \leftskip=14pt \hskip 0mm {
``As the brain forms images of an object --- such as a face, 
a melody, a toothache, the memory of an event --- and as the 
images of the object \textit{affect} the state of the organism, 
yet another level of brain structure creates a swift nonverbal account
of the events that are taking place in the varied brain regions 
activated as a consequence of the object-organism 
interaction. The mapping of the object-related 
consequences occurs in first-order neural maps 
representing the proto-self and object; the account of 
the \textit{causal relationshi}p between object and 
organism can only be captured in 
second-order neural maps.  \dots one might say 
that the swift, second-order nonverbal account 
narrates a story: \textit{that of the organism
caught in the  act of representing its
own changing 
state as it goes about 
representing something else}."}\vskip 0mm}

We think we know what is happening now? One needs only 
attempt a translation of our `understanding' into 
an engineering context to realise we still have some 
way to go. 

Other hypercomputational models are even more resistant 
to deconstruction. Roughly speaking, the more plausible the 
proposed instance of hypercomputational context, the 
more clearly present are all the ingredients we have 
drawn out above, and 
the less amenable it is to deconstructive inspection. 
For instance, we have Etesi and Nemeti's closely argued 
paper \cite{nemeti}, describing how relativistic 
considerations (involving the actuality of such  
things as large rotating black holes in 
galactic nuclei) may lead to effectively computable 
functions which are not 
Turing computable. The argument even puts up a  
credible resistance to Martin Davis'  destructive 
agenda in \cite{davis}, but finding 
out exactly what is hidden behind the complexities 
of the nonlinear mathematics involved is 
a hard task.\footnote{Although the new paper of 
N\'emeti and D\'avid  in the
present  volume  seeks to explain more thoroughly how and why such 
 general relativistic computers work.} This tends to be true, of course, 
of all the real-world, physically based examples, 
starting with the suggestion of Georg Kreisel 
\cite[p.143]{kreisel},  that a
collision problem related to
the 3-body problem might
give ``an analog computation of a
non-recursive function (by repeating collision
experiments sufficiently often)."

\section{Parallelism and interconnectivity}\label{interconnect}

So the building of routes through the Turing barrier appear blocked 
by parallel barriers to deconstruction of familiar processes. 
It is like exploration without maps and other navigational 
aids. Adventurous investigators may on occasion 
happen on strange and unrecognised computational 
territory, but when they return, their `discoveries' do not   
fit with the familiar scientific landscape, and there is no 
consensus about how to interpret their accounts. This certainly 
applies to how we view the role of parallelism and 
 interconnectivity 
within computational processes. 

As early as 1988  
Paul Smolensky  observed in his 
influential \textit{Behavioral and Brain Sciences\/} 
paper \cite[p.3]{smolensky} that: 

{\vskip 0mm\noindent\rightskip=14pt \leftskip=14pt \hskip 0mm {\sl
``There is a reasonable chance that connectionist 
models will lead to the development of new 
somewhat-general-purpose self-programming, 
massively parallel analog computers, and 
a new theory of analog parallel computation: 
they may possibly even challenge the 
strong construal of Church's Thesis as the 
claim that the class of well-defined computations 
is exhausted by those of Turing machines."\vskip 0mm}}

The classical computability theorist will be 
wary of elevating what may be to the applied scientist very 
important operational considerations into something more 
fundamental.   
One has to be very careful about claiming that some 
example of computational parallelism cannot be 
simulated by a Turing machine. As is well-known (see, for example, 
David Deutsch \cite[p.210]{deutsch}), 
the massive parallelism delivered by   
quantum computation  as it is currently abstracted 
is perfectly well contained within the classical 
sequential model. What happens when 
one  raises the  level of 
internal connectivity? 
It is still not clear 
that without analogue elements, such as in real world physical 
processes, you cannot come up with adequate descriptions of 
the algorithmic content of a given connectionist model 
of computation. Maybe external interaction will 
leave the classical model floundering? This idea 
runs through a number of hypercomputational 
proposals, including Copeland's \cite{copeland98} 
rediscovery of oracle Turing machines. 
But as Davis so effectively did, all our classical 
computability theorist has to do is to 
widen the definition of ``internal" --- that is to make  the modelling   
more inclusive  ---  
 to shepherd the 
proposed new paradigm back into the classical fold. 

There are limits to this though, determined by those 
on the success of hypercomputational modelling. 

Goldin and Wegner \cite{GW2} certainly have the 
modelling of more than internal 
 connectivity in mind when they  quote  from   
Robin Milner's  
1991 Turing Award lecture \cite[p.80]{milner}:

{\rightskip=14pt \leftskip=14pt \hskip 0mm {
``Through the seventies, I became convinced that a theory of 
concurrency and interaction requires a new conceptual framework, not just a
refinement of what we find natural for sequential  
computing."}\vskip
0mm}

In \cite{GW}, for instance, Goldin and Wegner are not just talking
about  parallelism. 
And the inherent 
vagueness of  examples  they quote  both stretch 
the mathematical analysis, and the reductionist agenda which feeds 
on that, to its limits: 

{\rightskip=14pt \leftskip=14pt \hskip 0mm {
``One example of a problem that is not algorithmic 
is the following instruction from a recipe \cite{knuth}: 
`toss lightly until the mixture is crumbly.' This problem 
is not algorithmic because it is impossible 
for a computer to know how long to mix: this may depend on 
conditions such as humidity that cannot be predicted with 
certainty ahead of time. In the function-based 
mathematical worldview, all inputs must 
be specified at the start of the computation, 
preventing the kind of feedback that would be necessary 
to determine when it's time to stop mixing."}\vskip
0mm}

Our intention here is not to rehash the familiar over-worked  
arguments over whether hypercomputation exists or not. 
In section \ref{intro} we settled on  reasons to bypass those.    
I just want to give us reasons to look more closely at how the real world 
expresses those  aspects 
of hypercomputational mental activity which we 
isolated in section \ref{neural}, and for which a 
better understanding 
may well be necessary for the next computing revolution. 

\section{Hypercomputation computationally simulated}

It is when we look more closely at the operational 
benefits of connectionist models that we get an inkling 
of what is new in a more basic sense. 

Viewing schematically, we do not have to 
look further than Turing machines to observe 
computably generated objects which are not themselves 
computable. Shine a bright light on the graph of a computable 
function over the natural numbers, and its shadow is likely to be 
incomputable --- or more precisely, the projection of a 
suitable computable 
binary relation over the numbers produces an incomputable set 
of numbers. So it seems we already knew that you get incomputable 
objects by selectively observing algorithmic processes. 
The only 
problem with our Turing machine example is one seems to 
need infinitely much time, which makes the 
computable simulation of incomputability interesting, 
but difficult to connect with our own world. 

However, if we put a little flesh on the algorithmic 
skeleton that is the Turing machine, we can get something 
which looks of more immediate interest. We can view fractals, 
and other mathematical objects generated according to simple 
algorithmic rules   as 
computer simulations. The 
 extraordinary richness of structure we observe is matched 
by the as yet unsolved problems of showing that 
aspects of, for instance, the Mandelbrot and certain Julia sets are 
computable (see  \cite{hertling},  
and, for contrast, \cite{RW}, \cite{braverman}, \cite{rettinger}). 
Of course, such mathematical examples provide a metaphor 
for the way real-world complexity is generated by 
the iteration of simple algorithmic content, and 
nicely illustrate how the analysis of apparently finitary objects 
rests on uncompleted infinities and the mathematics of the infinite. 

There is a fruitful two-way dialectic based on this 
interplay between computation and natural phenomenon. 
Nature is systematically raided nowadays for 
computational metaphors. Implicit in this is an 
assumption that nature can do something for us not 
yet delivered by unaided computational models. This 
`something' may be of a purely operational character, 
or may well rest on the sort of simulation of 
hypercomputation talked about above. Whatever it is, 
it is enough to fuel the huge current interest in 
areas such as artificial neural networks, cellular automata (derived from 
self-reproduction), 
membrane computing (based on biological cells), 
evolutionary computing, L-systems (multi-cellular 
organisms), swarm computing, molecular computing, and 
so on. One does not even need to extract the metaphor 
to computationally exploit this `something' --- one can 
ask the natural world itself to do part of the 
computation. Of course, the analogue approach 
tends to niche applications, and is unlikely to ever 
deliver a relacement for our desktop computer! 
In general, a lot of ingenuity needs to go into 
making natural phenomena compute more than themselves. 

Both the computational metaphors, and the direct 
exploitation of natural phenomena, seem to suffer from 
a basic deficiency, for different reasons. The full 
complexity of human problem solving depends on 
the breaking up of problems into parts which can be attacked 
with varying techniques, a sort of `divide and conquer' 
approach, involving loops in reasoning and 
interactions between conceptual frameworks. This lack  
of modularity and subroutines seems to be a feature of natural 
computation 
as we know it, particularly clearly in the case of quantum 
and relativistic 
computation,  or analogue computing.  It is not that 
these features are absent from the real world. 
Far from it. But because we do not fully understand 
how they arise, we cannot replicate them in our 
 models, nor ensure stability in direct exploitation 
of their physical embodiments. 

An important footnote to what we have said above is 
that  
one does not have to limit usable simulations to 
those classically modelled via enumerations of data. 
And if we think in terms of mental visualisations, this is just as 
well.  As Turing pointed out (see Hodges \cite[p.361]{hodges}) in a talk 
 to the London Mathematical Society, in February 20, 1947:

{\rightskip=14pt \leftskip=14pt \hskip 0mm {\sl 
``\dots if a machine is expected to be infallible, it 
cannot also be intelligent. There are several 
theorems which say 
almost exactly that."}\vskip 0mm}
 
It is reasonable to think in terms 
of appropriate  levels of adjustment and approximation 
for particular simulations, all capturable 
in classical computability theoretic hierarchies. 
This does
not just mean that factual mistakes may need 
correcting, but refers to the whole process of 
detailed imaging of objects external and internal 
to the organism.

\section{Definability and   
hypercomputational effects observed}

The picture we have built up seems to be a paradoxical 
one. On the one hand we believe that the classical 
dichotomy between computable and computably enumerable 
sets is reflected in the real world in the 
chasm between  
phenomena we can compute and those we can simulate.   
We believe that a simulation does hypercomputationally that a Turing 
machine  
does not --- produces observable infinity in finite time. 
But for this to have a real impact on our conception of what  
happens in  reality, we need valid models\footnote{By which we mean 
abstract, mathematical models, which do not appeal to undeconstructed  elements 
of the reality they are intended to model.} 
that actually compute 
relative to the non-local heightened information 
content. To be precise, we have models which 
elevate information content, and models which 
exploit this, but no overall cohesive picture. And it seems 
that the problem we had understanding the mechanics of 
consciousness and the local representation of non-local 
neural phenomena in the brain has a precise parallel 
in the difficulty we have knitting together these two types 
of computational model. If you come up with a model 
which elevates information content (such as in 
Beggs and Tucker \cite{BT}), you will probably be told 
by someone that it is unrealisable. If you admit  a rich 
pre-existent information content in your model 
to give scope to the interactivity we identified earlier  
 (as do  
most of Davis' targets in  \cite{davis04}),  the sceptic 
will have enough to deride its \textit{hyper}computational 
content as also being pre-existent. Let us give our 
paradox a name:
\smallskip

\textbf{Davis' Paradox:}  {Hypercomputation 
is dependant on interactivity within a domain of   
rich information content. But the existence of rich information 
content depends on hypercomputation.}
\smallskip

Returning to our section \ref{neural} neural avatar of hypercomputation, 
notice that we find it hard to locate intelligence as 
something that resides purely within the 
autonomous brain. The richness of mental activity appears as an 
extension of an equally complex world with which it interacts. 
Here is Brooks \cite{brooks2} again, drawing out this 
inclusiveness in a particulary dramatic way:

{\rightskip=14pt \leftskip=14pt \hskip 0mm {
``Real computational systems are not rational agents that take inputs, 
compute logically, and produce outputs \dots It is hard to draw the 
line at what is
intelligence and what is environmental interaction. In a sense, 
it does not really
matter which is which, as all intelligent systems must be 
situated in some world or
other if they are to be useful entities."}\vskip 0mm}

In fact, Brooks \cite[p.139]{brooks3} plausibly argues  that 
there is a realistic approach to AI involving no 
internally generated representations, but 
rather using ``the world as its own model".

Let us now find a way past Davis' paradox, in a way that 
validates the intuitions of both Martin Davis and  
Jack Copeland.  And gives us  
a mathematical model which translates    
the dichotomy we started this section with into one between 
nature captured by linear mathematics, and emergent 
phenomena (see \cite{cooper.emerg}). The idea is 
to provide a universe of scientific information content 
within which it makes sense to situate 
algorithmic content sufficient to model that 
extractable from the material world. If you like, 
we want to provide a playing field upon which the 
scientifically observed game of existence can be 
played out. This kind of model will be 
an extension of the classical model correctly favoured by Davis, and with 
the same kind of potential inclusiveness, 
but more in tune with what we actually observe. 
And this model, rather than being over-stretched by Goldin and Wegner's 
examples of real-world interactivity  
referred to in section \ref{interconnect},  will 
positively lap them up. 

Our model has as its underlying domain the set of all real 
numbers, within which the scientist commonly describes 
the material universe. The algorithmic content the 
scientist is intent on extracting is inclusively 
modelled in terms of Turing's extension of his 
machine model of computation to one capturing 
computation \textit{relative} to auxiliary 
information, accessed piecemeal via real numbers, 
and outputting real numbers again. The real 
numbers involved are dealt with via finitary approximations, 
of course. One can find the details of this  
extended Turing model --- what is often called 
the \textit{Turing universe\/}, familiar to us 
via its associated \textit{Turing degree structure\/} --- 
in most current computability theory texts (such as 
\cite{cooper04}, \cite{od89}, or \cite{soare}). 

The intention now is to use this  
mathematical model, in conjunction 
with our scientific understanding 
of the workings of the neural setting, to, on the one hand, 
give precision to interpretations of what is happening in the brain, 
and on the other, to reinforce our perceptions  
of the relevance of the Turing model, and to 
refine those of the model's  
inevitable limitations. There are three 
obvious qualifications to be noted in claiming 
relevance for the extended Turing model. The first and 
most fundamental qualification is that the 
model can only help us 
understand the underlying information content 
of a given computationally complex environment, and that 
one must look more carefully to see what this means 
in terms of constraints on the detailed workings of 
the organism or physical context in question. 
The second is that the information content present 
in a given physical example will at most present some 
relatively small substructure  of the full mathematical 
model. And perhaps most problematically, there will always 
be concern about what part of the algorithmic content of the 
Turing universe is relevant. Considerations of time, space and 
observed algorithmic content may entail refinements of the model, 
and even involve movement to different 
but related models which maintain the 
overall analytical framework, but lead to rather different 
detailed conclusions. But let us get back to 
our intended dialectic between mathematical and physical. 

Firstly, 
note how all our predictively implementable 
 understanding of the world is via
its  algorithmic content. From a post-Church-Turing  
perspective, this is almost
tautological.  And most of what we scientifically understand  
in a wider sense is built on that 
algorithmic content. Our higher order descriptions of reality 
tend to be shaped  
within this predictive framework.  
It may be that the weather, say, is not highly predictable in 
detail, but we do believe that climatic intricacies   
are emergent from a context whose local ingredients 
have well-understood algorithmic content. This is not to 
say there are not scientific puzzles, which have still to be brought 
within this dominant overall picture. Viewing the sub-atomic 
world from above, we only see the emergence in the form 
 constrained randomness, while the underlying 
local algorithmic content and its connections to what we observe 
cannot be easily grasped. But the neural case is more accessible 
(which is partly why we chose it), and offers a more immediately
productive environment in which to work out parallels between 
mathematical explanation and apprehensions of physical processes.

What is most interesting about the Turing universe 
are questions about exactly which global relations are 
definable in terms of its local structure --- given by  the locally applied
algorithmic content, or 
\textit{Turing reductions\/} --- using natural 
language familiar to us from its precise formulations 
in first-order logic\footnote{The underlying concept of 
definability, basic to this article,  is that familiar to logicians,   
where a relation on a structure is termed {\it definable\/} if it can 
be described in terms of the atomic relations of the structure 
using formal language, usually first-order, based on that of everyday
usage.}. What is particularly 
puzzling about the human mind is how it 
manipulates mental images, formed non-locally 
within the brain,  into higher-order 
conceptual formations. The use of language, peculiar to 
our own species, seems to be connected to this ability, 
but not a prerequisite for it. Now, on the 
one hand, the fact that language can be 
used to capture such higher-order formations 
points to a connection with language and definability. 
On the other, the well established  fact  
 that language is not actually 
necessary for this, makes us look more closely at what 
a description or definition actually provides 
in the way of structured understanding. The background to 
this will be our earlier discussion \cite{cooper.emerg} of the link between 
 emergence and definability, within which we 
seek to situate our discussion of the  gap 
in our understanding of neural 
representation.

The claim is that a description gives us a 
handle on a complex emergent relation, in that the 
form of the description provides a recipe for the 
simulation of the phenomenon described. The brain may 
store a number of images in the form of non-local 
simulations. It may also store descriptions of how to 
access those images. Those descriptions must be in some sense 
centrally located, and involve mental images of  
how to algorithmically access mental features which may be far from 
algorithmic in their origins and current simulative character. 
In doing this, the brain does not necessarily use language.  
What is  involved is yet another 
layer of cognitive function, with what we can 
think of as a recipe given by the description  
 imaged, and becoming part of what 
Damasio calls core consciousness.  
More graphically: What is happening is that the 
recipe embodied in the definition of the emergent 
relation which is the location of the mental image, and which 
provides the means
to reconstructing  the relation and its link to the core, makes the relation 
part of a framework within which it can be recalled and 
become part of a hierarchical interactive network of 
concepts, ideas, facts, and other stored mental images. 
What is striking is that when we examine how 
we ourselves construct and implement 
such recipes, we notice that the 
related image has a `feeling'. It is this 
which seems to identify, and to be the subjective 
apprehension of, its entry into a core neural
structure,  and provides the means to what   
Damasio identifies as our `ownership' over the 
emergent relation via its description. 
Here is how Damasio \cite[p.148]{damasio} puts it:

{\rightskip=14pt \leftskip=14pt \hskip 0mm {
``The perspective for a melody you hear or 
for an object you touch is, quite naturally, the 
perspective of your organism because it is drawn on the 
modifications that your organism undergoes 
during the events of hearing or touching. As for the sense 
of ownership of images and the sense of agency over those 
images, they, too, are a direct consequence of the 
machinations which create perspective. They are inherent 
in those machinations as foundational sensory 
evidence. Later, our creative and educated brains 
eventually clarify the evidence in the form of 
subsequent inferences, which also become known to us."}\vskip 0mm}

But where does this structure, which is 
accompanied by such symptomatic feelings, come from? And how  
does it relate to our mathematical model? 
Within the mind and its wider context of the physical 
universe, the core consciousness which is the key to 
this structuring is itself immanently emergent, 
while at the same time being the key to new  
activities. In a sense, 
the core, Damasio's `core consciousness', is itself 
defined by its own activities and the  
storing of descriptions of its relationship to 
neural formations (\cite[p.183]{damasio}):  

{\rightskip=14pt \leftskip=14pt \hskip 0mm {
``Beyond providing a feeling of 
knowing and enhancement of the object, the images 
of knowing, assisted by memory and reasoning, 
form the basis for simple nonverbal inferences 
which strengthen the process of core consciousness. 
These inferences reveal, for instance, the close linkage 
between the regulation of life and the 
processing of images which is implicit in the sense of 
individual perspective. Ownership is hidden, 
as it were, within the sense of perspective, ready to 
be made clear when the following inference can be made: 
if these images have the perspective of this body I now 
feel, then these images are in my body --- they 
are mine."}\vskip 0mm}

What the extended Turing model provides is the 
deconstruction of specific mental processes into what may be 
characterised as auxiliary (internally or externally 
originating) 
hypercomputational input, algorithmic content (in terms of 
locally specific neural connections), and the recording  
and utilisation of relevant descriptions representing 
definitions in terms of the more global neural environment 
of mental inter-relationships.

\section{Conclusions}

The closest anyone  has got so far to 
actual computers with  
recognisably hypercomputational ingredients is by 
surfing physical reality in some way. There is 
a widespread suspicion that the world 
cannot be satisfactorily located within the standard 
computational model. This is not necessary for one to have an 
interest in new computational paradigms --- as in the case of 
quantum computation, there may be very important 
operational benefits, even though 
there is an underlying classical model. But this suspicion 
gets stronger the more difficult it is to divorce ones 
computational approach from its real-world origins. 

One way forward is to utilise the physical world's  
rich potential for computation, without worrying too 
much about understanding the underlying rules of the 
game. The likely success of this approach may be limited --- 
it takes ingenuity\footnote{Or, as in the case of relativistic 
computation,  extreme demands on 
resources, as surmised by N\'emeti and  D\'avid \cite{ND}.}  
 to get a natural
process to compute more  than itself ---  
but may be the best we can do in the short to medium term. Or ever! 

Let us take an analogy. The domestication of horses 
around five or six thousand years ago brought a 
revolution in transportation, only 
achieved through a creative interaction between humans and 
the  natural world. At that time, trying to understand the 
principles underlying the equine organism in order to 
synthesise an artificial horse was unthinkable. 
But a few thousand  years later   
there was enough understanding of  scientific basics 
to underpin the invention of the ``iron horse", leading, amongst 
other things, to 
the opening up of many previously isolated parts of the world to 
people with no riding skills whatsoever. 

The moral of this seems to be that perhaps there is a great 
deal to be got from an ad hoc computational relationship 
with the real world. But that we should not be daunted by  the  
sheer wonder natural structures inspire in us. It may be that 
the human brain, as an emergent phenomenon, has an intimate 
relationship with processes which are not easily simulable 
over significantly shorter time-scales than those 
to which natural evolution is subject --- a kind of 
earthy victory of the functionalist over the 
computationalist view of the human mind. Maybe we will 
never build an artificial brain, anymore than we can 
make an artificial horse. But this does not mean we may not one day 
have a good enough understanding of basic 
hypercomputational principles to build computers --- 
or firstly non-classical mathematical models of 
computation --- which do things undreamt of 
today. 

Rodney Brooks \cite[p.139]{brooks3} gets to have the last word:

{\rightskip=14pt \leftskip=14pt \hskip 0mm {
``I, and others, believe that human level
intelligence is too complex and little understood to be
correctly decomposed into the right subpieces at the
moment and that even if we knew the subpieces we
still wouldn't know the right interfaces between
them. Furthermore, we will never understand how to
decompose human level intelligence until we've had a
lot of practice with simpler level intelligences."}

\end{document}